\documentclass[12 pt]{amsart}
\usepackage{a4}
\usepackage{amssymb}

\newtheorem{thm}{Theorem}[section]
\newtheorem{lem}[thm]{Lemma}
\newtheorem{prop}[thm]{Proposition}
\newtheorem{cor}[thm]{Corollary}
\theoremstyle{definition}

\newtheorem{example}[thm]{Example}
\theoremstyle{remark}

\numberwithin{equation}{section}

\begin{document}

\title[Boundaries for algebras of holomorphic functions]{Boundaries for algebras of holomorphic functions on Banach spaces}

\author{Yun Sung Choi}
\address{Department of Mathematics, POSTECH, San 31, Hyoja-dong,
Nam-gu, Pohang-shi, Kyungbuk, Republic of Korea, +82-054-279-2712}

\email{mathchoi@postech.ac.kr}

\thanks{This work was supported by grant  No. R01-2004-000-10055-0 from the
Basic Research Program of the Korea Science \& Engineering
Foundation.}
\thanks{The third named author is the corresponding author.}

\author{Kwang Hee Han}
\email{hankh@postech.ac.kr}

\author{Han Ju Lee}
\email{hahnju@postech.ac.kr}





\subjclass[2000]{46E50, 46B20, 46B45}


\keywords{boundary for algebra, Shilov boundary,  complex convexity,
local uniform monotonicity, Banach sequence space.}

\begin{abstract}
We study the relations between boundaries for algebras of
holomorphic functions on Banach spaces and complex convexity of
their balls. In addition, we show that the Shilov boundary for
algebras of holomorphic functions on  an order continuous sequence
space $X$ is the unit sphere $S_X$ if $X$ is locally c-convex. In
particular, it is shown that the unit sphere of the Orlicz-Lorentz
sequence space $\lambda_{\varphi, w}$ is the Shilov boundary for
algebras of holomorphic functions on $\lambda_{\varphi, w}$ if
$\varphi$ satisfies the $\delta_2$-condition.
\end{abstract}

\maketitle

\section{Introduction and preliminaries}
Let $X$ be a complex Banach space and let $B_X$ be the closed unit
ball of $X$. We  denote  by $H(B_X)$ the set of all holomorphic
functions on the interior of $B_X$, and by $C_b(B_X)$ the Banach
algebra of bounded continuous functions on $B_X$ with the sup norm.

Globevnik \cite{G1} defined and studied the following analogues of
the classical disc algebra:
\begin{align*} \mathcal{A}_b(B_X) &=\left\{ f\in H(B_X): f \in C_b(B_X)  \right\}
\\
\mathcal{A}_u(B_X) &= \left\{ f\in \mathcal{A}_b(B_X) : f \mbox{ is
uniformly continuous on } B_X \right\}\end{align*}

It is shown in \cite{ACG} that $\mathcal{A}_u(B_X)$ is a proper
subset of $\mathcal{A}_b(B_X)$ if and only if $X$ is an infinite
dimensional Banach space. Then it is easy to see that both
$\mathcal{A}_b(B_X)$ and $\mathcal{A}_u(B_X)$ are Banach algebras
when given the natural norm
\[ \|f\| = \sup\{ |f(x)| : x\in B_X\}.\]

Let $K$ be a Hausdorff topological space and $\mathcal{A}$ a {\it
closed function algebra} on $K$, that is, a closed subalgebra of
$C_b(K)$. A subset $F$ of $K$ is called a {\it boundary} for
$\mathcal{A}$ if for all $f\in \mathcal{A}$ we have
\[ \|f\| = \sup_{x\in F} |f(x)|.\] If the intersection of all
closed boundaries for $\mathcal{A}$ is again a boundary for
$\mathcal{A}$, then it is called the {\it Shilov boundary} for
$\mathcal{A}$, denoted by $\partial \mathcal{A}$. We recall that a
function algebra $\mathcal{A}$ is said to be {\it separating} if (i)
for two distinct points $x, y$ in $K$, there is an element $f\in
\mathcal{A}$ such that $f(x) \neq f(y)$ and (ii) for each $t\in K$
there is a $f\in \mathcal{A}$ such that $f(t)\neq 0$. A {\it uniform
algebra} on a compact Hausdorff space $K$ is a closed function
algebra
 which contains constants and separates the points of $K$.

Given a closed function algebra $\mathcal{A}$ on a metric space $K$,
a set $S\subset K$ is called a {\it peak set} for $\mathcal{A}$ if
there exists $f\in \mathcal{A}$ such that $f(S)=1,~~ |f(x)|<1
~~(x\in K\setminus S)$. A set $S\subset K$ is called a {\it strong
peak set} for $\mathcal{A}$ if there exists $f\in \mathcal{A}$ such
that $f(S)=1$ and for every $\epsilon >0$ there exists $\delta >0$
with $|f(x)|<1-\delta$ whenever ${\rm dist}(x, S)>\epsilon.$ If $S$
consists of only one point $p$ and if it is a {\it peak set} (resp.
{\it strong peak set}) for $\mathcal{A}$, then the point $p$ is
called a {\it peak point} (resp. {\it strong peak point}) for
$\mathcal{A}$.

The set of all peak points for $\mathcal{A}$ is called the {\it
Bishop boundary} for $\mathcal{A}$ and denoted by $\rho
\mathcal{A}$. Note that if $K$ is compact, then a peak point $x\in
K$ for $\mathcal{A}$ is also a strong peak point for $\mathcal{A}$,
hence every closed boundary for $\mathcal{A}$ contains $\rho
\mathcal{A}$. Further, Bishop \cite[Theorem 1]{B} showed that for
any uniform algebra $\mathcal{A}$ on a compact metrizable space $K$,
\[\partial \mathcal{A} = \overline{\rho \mathcal{A}}.\]

Given a convex set $M\subset X$,  a point $x\in M$ is called a {\it
real} (resp. {\it complex}) {\it extreme point} of $M$ if for every
nonzero $y\in X$, there is a real (resp. complex) number $\zeta$
such that $|\zeta|\le 1$ and $x+\zeta y \notin M$. The set of all
real (resp. complex) extreme points of $M$ is denoted by
$Ext_{\mathbb{R}}(M)$ (resp. $Ext_{\mathbb{C}}(M)$). Let
$\mathcal{A}$ be a uniform algebra on a compact Hausdorff space $K$.
Let $\mathcal{A}^*$ be the dual Banach space of $\mathcal{A}$ and
let $S_1^*$ be the intersection of the unit sphere
$S_{\mathcal{A}^*}$ of $\mathcal{A}^*$ with the hyperplane $\{
x^*\in \mathcal{A}^* : x^*(1) = 1\}$. The set $\chi \mathcal{A} = \{
x\in K: \delta_x \in Ext_{\mathbb{R}}(S_1^*)\}$ is called the {\it
Choquet boundary} for $\mathcal{A}$.

It is well-known (see \cite[Theorem~9.7.2]{L}) that if $\mathcal{A}$
is a uniform algebra on a compact metrizable space $K$, then
\begin{equation} \label{eqz11}\rho \mathcal{A} = \chi \mathcal{A}.\end{equation}

Given a convex compact subset $K$ in a complex locally convex space
$E$, Arenson \cite{A} considered the uniform algebra
$\mathcal{P}(K)$ generated by the constants and restrictions to $K$
of functions from $E^*$, and showed that
\begin{equation} \label{eqz12}\chi \mathcal{P}(K) =
Ext_{\mathbb{C}}(K).\end{equation} In particular, if $K$ is
metrizable, we have
\begin{equation} \label{eqz13}\rho \mathcal{P}(K) = \chi
\mathcal{P}(K) = Ext_{\mathbb{C}}(K).\end{equation}

On the other hand, it is shown in \cite{G1} that $\rho
\mathcal{A}_b(B_X)\subset Ext_\mathbb{C}(B_X)$. We also note that
every closed boundary for a function algebra $\mathcal{A}$ must
contain the set of all strong peak points for $\mathcal{A}$.

When $X$ is finite dimensional, we get the following observation
from the above and $\rho \mathcal{P}(B_X)\subset \rho
\mathcal{A}_b(B_X).$

\begin{prop}\label{finite}
If $X$ is finite dimensional, then $$\rho \mathcal{P}(B_X) = \chi
\mathcal{P}(B_X) = Ext_{\mathbb{C}}(B_X)=\rho \mathcal{A}_b(B_X),$$
and $$
\partial \mathcal{A}_b(B_X) = \overline{Ext_\mathbb{C}(B_X)}.$$
\end{prop}

 The following two observations are easy.

\begin{prop}\label{propinfinitedense}
Let $X$ be a Banach space and suppose that the set of all strong
peak points for $\mathcal{A}_u(B_X)$ is dense in the unit sphere
$S_X$. Then the unit sphere is the Shilov boundary for both
$\mathcal{A}_u(B_X)$ and $\mathcal{A}_b(B_X)$.
\end{prop}

\begin{prop}
Let $X$ be a Banach space and suppose that the unit sphere $S_X$ is
the Shilov boundary for $\mathcal{A}_u(B_X)$. Then a subset $F$ of
$B_X$ is a boundary for $\mathcal{A}_u(B_X)$ if and only if it is a
boundary for $\mathcal{A}_b(B_X)$.
\end{prop}

For each $x\in B_X$, we call the set $\mathcal{F}(x)$ the {\it face}
at $x$, which is defined by
\[ \mathcal{F}(x) = \left\{x+y : y \in X,\ \sup_{|\zeta| \le 1} \| x+ \zeta y\| \le 1\right\}.\]

Notice that $x\in S_X$ is a complex extreme point of $B_X$ if and
only if $\mathcal{F}(x)=\{x\}$. A Banach space $X$ is said to be
{\it strictly c-convex} if every point of $S_X$ is a complex extreme
point.

By the maximum modulus theorem, we obtain the following
\begin{prop}\label{finitface}
Let $S$ be a peak set for $\mathcal{A}_b(B_X)$. Then for each $x\in
S$, $\mathcal{F}(x)$ is contained in $S$.
\end{prop}

This shows that every peak point for $\mathcal{A}_b(B_X)$ is a
complex extreme point of $B_X$, which is Theorem~4 in \cite{G2}.

A point $x\in S_X$ is said to have a {\it strong face} if for each
$\epsilon>0$, there is $\delta(\epsilon)>0$ such that if ${\rm
dist}(\mathcal{F}(x), y)\ge \epsilon$, then
\[ \sup_{0 \le \theta\le 2\pi} \|x+e^{i\theta}(x-y)\| \ge 1 +
\delta(\epsilon).\] A Banach space $X$ is said to be {\it locally
c-convex} if it is strictly c-convex and every point of the unit
sphere $S_X$ has a strong face. The maximum modulus theorem shows
that if two elements $x, y$ in a Banach space satisfy $\sup_{0\le
\theta\le 2\pi} \|x+e^{i\theta}y\|\le M$, then $\max\{\|x\|,
\|y\|\}\le M$.

\begin{prop}
Suppose that $X$ is a finite dimensional Banach space. Then every
point of $S_X$ has a strong face.
\end{prop}

\begin{proof}Suppose otherwise. Then there exist $x\in S_X$,
$\epsilon>0$ and a sequence $\{y_n\}$ in $X$ such that ${\rm
dist}(\mathcal{F}(x), y_n)\ge \epsilon$ for each $n$, but\[ \lim_{n
\rightarrow \infty} \sup_{0\le \theta\le 2\pi}
\|x+e^{i\theta}(x-y_n)\| =1.\] So we get
\[ M := \sup_{n \in \mathbb{N}} \sup_{0\le \theta\le 2\pi}
\|x+e^{i\theta}(x-y_n)\| < \infty.\]
 Hence $\sup_n
\|x-y_n\|\le M$. So we may assume that $y_n$ converges to $y$. Then
${\rm dist}(\mathcal{F}(x), y)\ge \epsilon$. For each $\theta\in
\mathbb{R}$,
\[ \|x+ e^{i\theta}(x-y) \| \le \|x+ e^{i\theta}(x-y_n)\| + \|y_n -
y\|.\] This shows that
\[ \sup_{0\le \theta\le 2\pi} \|x+ e^{i\theta}(x-y) \| \le
\lim_n\sup_{0\le \theta\le 2\pi} \|x+ e^{i\theta}(x-y_n) \| =1.\]
Therefore, $y = x+ (y-x)\in \mathcal{F}(x)$, which contradicts ${\rm
dist}(\mathcal{F}(x), y)\ge \epsilon$.
\end{proof}

The modulus of complex convexity of a complex Banach space $X$ is
defined by
\[ H_X(\epsilon) = \inf \left\{ \sup_{0\le \theta\le 2\pi} \|x+e^{i\theta} y\| -1
: x\in S_X, \|y \|\ge \epsilon      \right\}\] for each
$\epsilon>0$. A complex Banach space $X$ is said to be {\it
uniformly c-convex} if $H_X(\epsilon)>0$ for all $\epsilon>0$.  If
$X$ is uniformly c-convex, then every point in $S_X$ has a strong
face. A finite dimensional strictly c-convex space is uniformly
c-convex.

A sequence  $x=\{x(k)\}$ is said to be {\it positive} if $x(k)\ge 0$
for each $k\in \mathbb{N}$. We define a partial order $x\ge y$ if
$x-y$ is positive. The absolute value of $x$ is defined to be
$|x|=\{|x(k)|\}$. A {\it Banach sequence space} $(X, \|\cdot\|)$ is
a Banach space consisting of sequences satisfying the following : if
$x$ is a sequence with $|x|\le |y|$ for some $y\in X$, then $x\in X$
and $\|x\|\le \|y\|$. A Banach sequence space is said to be {\it
order continuous} if any sequence $\{x_n\}$ in $X$ satisfying
\[ 0\le x_1 \le x_2 \le \cdots \le y \ \ \mbox{for some positive}\ \ y\in X,\]
 is norm-convergent. The vector $e_j$ is defined to
have 1 in  the $j$-th component with all zeros in the other
components. Note that if a Banach sequence space $X$ is order
continuous, then $\{e_n\}$ is a basis of $X.$

It is known in \cite{Lee, Lee2} that a uniformly c-convex sequence
space is order continuous.

A Banach sequence space $X$ is said to be {\it strictly monotone} if
for every  pair $y\ge x\ge 0$ with $y\neq x$, we have $\|y\|>\|x\|$.
Recall also that a Banach sequence space $X$ is said to be {\it
lower (resp. upper) locally uniformly monotone} if for any positive
$x\in S_X$ and any $0<\epsilon<1$ (resp. $\epsilon>0$) there is
$\delta=\delta(\epsilon, x)>0$ such that the condition $0\le y \le
x$ (resp. $y\ge 0$) and $\|y\|\ge \epsilon$ implies
\[ \| x - y \| < 1-\delta \ \ \ \ \mbox{(resp.} \ \  \|x+ y\|\ge 1
+\delta).\] A Banach sequence space $X$ is said to be {\it uniformly
monotone} if given $\epsilon>0$, there  is $\delta(\epsilon)>0$ such
that
\[ \inf\{ \|\ |x| + |y| \ \| : \|y\|\ge \epsilon, x\in S_X\} \ge 1 +
\delta.\] A uniformly monotone Banach sequence space is both lower
and upper locally uniformly monotone.

It is shown in \cite[Theorem 1]{FK} that a Banach sequence space is
lower locally uniformly monotone if and only if it is strictly
monotone and order continuous. It is also shown in \cite{HN, Lee,
Lee2} that a Banach sequence space is strictly (resp. uniformly)
monotone if and only if it is  strictly (resp. uniformly) c-convex.

\section{Boundaries of $\mathcal{A}_u(B_X)$ and $\mathcal{A}_b(B_X)$}

\begin{prop} Suppose that $X$ is a complex Banach space.
Let $F$ be a boundary for $\mathcal{A}_u(B_X)$ and let $P$ be a
norm-one projection with a finite dimensional range $Y$. Then
\[ Ext_{\mathbb{C}} (B_Y) \subset \overline{P(F)}. \]
\end{prop}

\begin{proof}
Suppose  $x_0 \in Ext_{\mathbb{C}}(B_Y)\backslash \overline{P(F)}$.
Then there exists $\epsilon_0 > 0$ such that $\|P(x) - x_0 \| \ge
\epsilon_0$ for every $x\in F$. By Proposition~ \ref{finite} $x_0$
is a strong peak point for the algebra $\mathcal{A}_u(B_Y)$, that
is, there is a $g\in \mathcal{A}_u(B_Y)$ such that $g(x_0)=1$ and to
every $\epsilon>0$ corresponds a $\delta(\epsilon)>0$ satisfying
\[ |g(y)| < 1-\delta(\epsilon),\]
for all $y\in B_Y$ with $\|y-x_0\|\ge \epsilon.$ Take $f= g\circ
P\in \mathcal{A}_u(B_X)$. Then $f(x_0)=1$ and for every $x\in F$ we
have
\[ |f(x)| = |g(P(x))| < 1-\delta(\epsilon_0).\] This
contradicts the fact that $F$ is a boundary of $\mathcal{A}_u(B_X)$.
\end{proof}

\begin{prop}\label{finitedimension}
Suppose that $X$ is a complex Banach space with the following
properties: There is a collection $\{P_{\alpha}\}_{\alpha\in A}$ of
projections $P_{\alpha}$ with finite dimensional ranges $Y_{\alpha}$
such that $\cup_{\alpha\in A} Y_{\alpha}$ is dense in $X$, and for
each $\alpha\in A$ \[\sup_{0\le \theta\le 2\pi}\|P_{\alpha} +
e^{i\theta} (I-P_{\alpha})\|\le 1.\]  Then a set $F\subset B_X$ is a
boundary for $\mathcal{A}_u(B_X)$ if
$$Ext_{\mathbb{C}}(B_{Y_{\alpha}})\subset
\overline{P_{\alpha}(F)}$$ for every $\alpha\in A$.
\end{prop}

\begin{proof}
 Suppose that $F$ is not a boundary for
$\mathcal{A}_u(B_X)$. Then there are $f\in \mathcal{A}_u(B_X)$,
$\|f\|=1$ and $\epsilon>0$ such that $|f(x)|<1-\epsilon$ for every
$x\in F$. $\rho \mathcal{A}_u(B_{Y})$ is a boundary for
$\mathcal{A}_u(B_{Y})$ if $Y$ is finite dimensional. Since
$\cup_{\alpha\in A} B_{Y_{\alpha}}$ is dense in $B_X$, it follows
from Proposition~ \ref{finite} that $$ \|f\| = \sup_{x\in
\cup_{\alpha} S_{Y_{\alpha}}} |f(x)| = \sup_{x\in \cup_{\alpha} \rho
\mathcal{A}_u(B_{Y_{\alpha}})} |f(x)| =\sup_{x\in \cup_{\alpha}
Ext_\mathbb{C}(B_{Y_{\alpha}})} |f(x)|.$$ Hence there is a sequence
$\{x_n\}$ such that \[x_n \in \bigcup_{\alpha} Ext_\mathbb{C}
(B_{Y_{\alpha}}),\ \ \ \mbox{and} \ \ \ \ \lim_n |f(x_n)|=1.\]

Because $f$ is continuous and $Ext_\mathbb{C}(B_{Y_{\alpha}})\subset
\overline{P_{\alpha}(F)}$ for every $\alpha\in A$, there is a
sequence $\{u_n\}$ such that
\[ \{u_n\} \subset \bigcup_{\alpha} P_{\alpha}(F), \ \ \ \ \mbox{and}  \ \ \ \
\lim_n |f(u_n)|=1.\] Each $u_n$ has the form $u_n = P_{\alpha_n}z_n$
where $z_n\in F$. Set $v_n= (I-P_{\alpha_n})z_n$, and $z_n = u_n +
v_n$. By the uniform continuity of $f$ there exists $\delta$,
$0<\delta <1$ such that if $\|x_1 - x_2\|_X \le \delta$ and $x_1,
x_2\in B_X$, then $|f(x_1) - f(x_2)|< \epsilon/2$. Thus we get for
every $n\in \mathbb{N}$,
\[ |f(u_n + (1-\delta)v_n) - f(u_n + v_n)| < \frac{\epsilon}2.\]

Further, since $z_n=u_n + v_n\in F$, we have $|f(u_n +
v_n)|<1-\epsilon$, and consequently for each $n\in \mathbb{N}$,
\begin{equation}\label{eq:contra} |f(u_n + (1-\delta )v_n)| < 1-
\frac{\epsilon}2.\end{equation}  On the other hand, since
$$\|P_{\alpha_n} + e^{i\theta} (I-P_{\alpha_n})\|\le 1$$ for every
$\theta\in \mathbb{R}$, the maximum modulus theorem shows that for
every $\zeta\in \mathbb{C}$ with $|\zeta|\le 1$,
\[ \left\|u_n + \frac{1}{1-\delta}\zeta [(1-\delta)v_n] \right\| \le
1.\] By  \cite[Lemma~1.4]{G1}, there is $C(\epsilon)<\infty$ such
that for  each $n\in \mathbb{N}$,
\[ |f(u_n + (1-\delta)v_n) - f(u_n)| <
C(\epsilon)(1-|f(u_n)|).\]

Since $\lim_n |f(u_n)|=1$, it follows that $\lim_n |f(u_n +
(1-\delta)v_n)|=1$, which contradicts (\ref{eq:contra}).
\end{proof}

\begin{cor} Suppose that $X$ is a complex Banach space with a
sequence $\{P_n\}$ of projections with the same properties as in
Proposition \ref{finitedimension}. Then a set $F\subset B_X$ is a
boundary for $\mathcal{A}_u(B_X)$ if and only if the closure of
$P_n(F)$ contains $Ext_{\mathbb{C}}(B_{Y_n})$ for every positive
integer $n$.
\end{cor}

We remark that an order continuous Banach sequence space has the
properties outlined in Proposition \ref{finitedimension}.

\begin{cor}\label{specialcase}
Let $X$ be an order continuous Banach sequence space. Let $F\subset
S_X$ and let $P_n$ be a sequence of coordinate projections with
finite dimensional range $Y_n$ such that every finite subset of
$\mathbb{N}$ is contained in the support of some $P_n$. If
$Ext_{\mathbb{C}}(B_{Y_n})\subset \overline{P_n(F)}$ for each $n\in
\mathbb{N}$, then $F$ is a boundary for $\mathcal{A}_u(B_X)$.
\end{cor}

Corollary~\ref{specialcase} extends Theorem~1.5 in \cite{G1}.

\begin{prop}\label{propbboundary}
Suppose that there is a family $\{Y_\alpha\}_{\alpha\in A}$ of
finite dimensional subspaces of a Banach space $X$ such that
$\cup_\alpha B_{Y_\alpha}$ is dense in $B_X$. Then the set
$F=\cup_\alpha Ext_{\mathbb{C}}(B_{Y_\alpha})$ is a boundary for
$\mathcal{A}_b(B_X)$.
\end{prop}

\section{Shilov Boundary for $\mathcal{A}_u(B_X)$  and $\mathcal{A}_b(B_X)$}

\begin{prop}\label{strongface}
Suppose that $x_0\in S_X$ has a strong face and $T$ is a bounded
operator of $X$ into $X$ with $Tx_0 =x_0$ and \[\sup_{0\le\theta\le
2\pi}\|T+e^{i \theta}(I-T)\|\le 1.\] Then for each $\epsilon>0$,
there is $\delta(\epsilon)>0$ such that whenever ${\rm dist}(F(x_0),
y)\ge\epsilon$ and $y\in B_X$, we get
\[ \|x_0 - Ty \|\ge \delta(\epsilon).\]
\end{prop}
\begin{proof}
Suppose, on the contrary, that there is $\epsilon_0>0$ such that
\[ \inf\{ \|x_0 - Ty\|: {\rm dist}(F(x_0), y)\ge\epsilon_0, \ y\in B_X\} =0. \] Then there
is a sequence $\{s_n\}$ in $B_X$ such that ${\rm dist}(F(x_0),
s_n)\ge \epsilon_0$ and
\[\lim_{n\rightarrow \infty} Ts_n = x_0.\]
Since $x_0$ has a strong face, there are  $\delta_1>0$ and
$\theta_n\in \mathbb{R}$ such that for every $n$, we have
\[\|x_0+e^{i\theta_n}(s_n-x_0)\| \ge 1 + \delta_1.\] So
\begin{align*} 1+\delta_1 &\le \|x_0 + e^{i\theta_n}(s_n -x_0)\|
\\&\le \|Ts_n +
e^{i\theta_n}(s_n-x_0)\| + \|x_0 - Ts_n\|\\
&\le \|Ts_n + e^{i\theta_n}(I-T)s_n)\| + 2\|x_0 - Ts_n\|.
\end{align*}  This implies that
\[ 1+\delta_1 \le \limsup_{n\rightarrow\infty} \|Ts_n + e^{i\theta_n}(I-T)s_n)\| \le
1,\] which is a contradiction.\end{proof}

\begin{prop}\label{peakset}
 Let $X$ be a complex Banach space and let $P$ be a projection of
$X$ onto a finite dimensional subspace $Y$ such that
\[\sup_{0\le\theta\le2\pi}\|P + e^{i\theta} (I-P)\|\le 1.\]
If $x_0\in Ext_\mathbb{C}(B_Y)$ and  $x_0$ has a strong face in $X$,
then $\mathcal{F}(x_0)$ is a strong peak set for
$\mathcal{A}_u(B_X)$.
\end{prop}

\begin{proof}
By Proposition~ \ref{finite}, $x_0$ is a strong peak point for
$\mathcal{A}_b(B_Y)$ since $Y$ is finite dimensional. Hence there is
a peak function $g\in \mathcal{A}_b(B_Y)$ such that $g(x_0)=1$ and
for each $\epsilon>0$ there is $\tilde{\delta}(\epsilon)>0$ such
that for every $\|x_0 - y\|\ge \epsilon$ and $y\in B_Y$, we have
\[ |g(y)| < 1-\tilde{\delta}(\epsilon).\]
By Proposition~\ref{strongface} we get $\delta(\epsilon)>0$ such
that if ${\rm dist}(\mathcal{F}(x_0),y) \ge \epsilon$
\[ \|x_0 - Py \| \ge \delta(\epsilon).\]
Take $f = g\circ P$. Then $f\in \mathcal{A}_u(B_X)$ and for each
$z\in B_X$ with ${\rm dist}(\mathcal{F}(x_0),  z)\ge \epsilon$, we
have $f(\mathcal{F}(x_0))=1$ by the maximum modulus theorem and $
|f(z)| <1- \tilde{\delta}( \delta(\epsilon))$. This implies that
$\mathcal{F}(x_0)$ is a strong peak set for $\mathcal{A}_u(B_X)$.
\end{proof}

Every element $x$ in the torus in $c_0$ has a strong face
$\mathcal{F}(x)$ and hence $\mathcal{F}(x)$ is a strong peak set for
$\mathcal{A}_u(B_X)$. The following generalizes Theorem~1.9 in
\cite{G1}.

\begin{prop}\label{propboundarypeakset}
Let $X$ be a complex Banach space as in Proposition
\ref{finitedimension}. Suppose that every point of $\cup_{\alpha}
Ext_\mathbb{C}(B_{Y_{\alpha}})$ has a strong face in $X$. Then
$F\subset B_X$ is a boundary for $\mathcal{A}_u(B_X)$ if and only if
${\rm dist}(F, S)=0$ for each strong peak set $S$ for
$\mathcal{A}_u(B_X)$.
\end{prop}

\begin{proof}
The necessity is clear. Conversely, suppose that there is a subset
$F\subset B_X$ such that ${\rm dist}(F, S)=0$ for each strong peak
set $S$ for $\mathcal{A}_u(B_X)$. We shall show that for each
$\alpha$ the closure of $P_{\alpha}(F)$ contains
$Ext_\mathbb{C}(B_{Y_{\alpha}})$. By applying
Proposition~\ref{finitedimension}, we get the desired result.

Now, let $x_0\in Ext_\mathbb{C}(B_{Y_{\alpha}})$. By
Proposition~\ref{peakset}, its face $\mathcal{F}(x_0)$ is a strong
peak set for $A_u(B_X)$. Hence ${\rm dist}(F, \mathcal{F}(x_0))=0$.
Then there are sequences $\{x_0 + y_k\}_k$ in $\mathcal{F}(x_0)$ and
$\{z_k\}_k$ in $F$ such that
\begin{equation}\label{eqlimitext}
\lim_{k\rightarrow \infty} \| (x_0 + y_k) - z_k\| =0.\end{equation}
Since $x_0 + y_k$ is in $\mathcal{F}(x_0)$ and $\|P_{\alpha}\| = 1$,
we have for each real $\theta$,
\[ \|x_0 + e^{i\theta} P_{\alpha}(y_k) \| \le \|x_0 + e^{i\theta} y_k\|\le 1.\]
Since $x_o$ is a complex extreme point of $B_{Y_{\alpha}}$,
$P_{\alpha}(y_k)=0$, and so (\ref{eqlimitext}) shows that
\[ \limsup_{k\rightarrow \infty}\| x_0 - P_{\alpha}(z_k)\| \le
\limsup_{k\rightarrow \infty} \|(x_0 + y_k) - z_k\| = 0.\] Therefore
$x_0$ is in the closure of $P_{\alpha}(F)$.
\end{proof}

In the proof of Proposition~\ref{propboundarypeakset}, it is
sufficient for $F\subset B_X$ to be a boundary for
$\mathcal{A}_u(B_X)$ that ${\rm dist}(F, \mathcal{F}(x))=0$ holds
for every $x\in \cup_{\alpha} Ext_\mathbb{C}(B_{Y_{\alpha}})$.

\begin{cor}\label{propz1001}
Let $X$ be a locally c-convex sequence space. Suppose that $x_0\in
S_X$ is finitely supported. Then $x_0$ is a strong peak point for
$\mathcal{A}_u(B_X)$. In particular, if in addition, $X$ is order
continuous, then the set of all strong peak points for
$\mathcal{A}_u(B_X)$ is dense in $S_X$.
\end{cor}
\begin{proof}
Suppose that $Y={\rm span}\{e_1, \ldots, e_n\}$ contains $x_0$.
Hence $x_0$ is a complex extreme point of $B_Y$. Let $P: X
\rightarrow Y$ be the projection defined by
\[ P(x(1), x(2) ,\ldots ) = (x(1), x(2), \ldots, x(n), 0, 0,\ldots
).\] Clearly $\|P + e^{i\theta}(I-P)\|\le 1$ for all $\theta\in
\mathbb{R}$. By Proposition~\ref{peakset}, $x_0$ is a strong peak
point for $\mathcal{A}_u(B_X)$.  Notice that if a Banach sequence
space is order continuous, then the set of all finitely supported
elements in $X$  is dense in $X$.\end{proof}

Proposition~\ref{propinfinitedense} and Corollary~\ref{propz1001}
show the following theorem.
\begin{thm}\label{thmsilov}
Let $X$ be an order continuous locally c-convex Banach space. Then
$S_X$ is the Shilov boundary for both $\mathcal{A}_u(B_X)$ and
$\mathcal{A}_b(B_X)$.
\end{thm}

By \cite{Lee, Lee2}, every uniformly c-convex sequence space is
order continuous.

\begin{prop} \label{propmonotonecomplex}
 A Banach sequence space $X$ is upper locally uniformly monotone if
and only if it is locally c-convex.
\end{prop}
\begin{proof}
 Suppose  $X$ is locally c-convex. Then for
each positive $x\in S_X$ and $\epsilon>0$ there is
$\delta=\delta(x,\epsilon)>0$ such that for all $z \in X$ with
$\|z\|\ge \epsilon$
\[ \sup_{0 \le \theta \le 2\pi} \|x+ e^{i\theta} z\| \ge 1 + \delta.\]
Hence we have for every $y\ge 0$ with $\|y\|\ge \epsilon$,
\[  \|x+ y\| \ge \inf \left\{\sup_{0 \le \theta \le 2\pi} \|x+ e^{i\theta} z\| : \|z\| \ge \epsilon \right\} \ge 1 +
\delta.\] So $X$ is upper locally uniformly monotone.

Conversely, suppose that $X$ is upper locally uniformly monotone. If
$x, y\in X$, then by \cite[Theorem~7.1]{DGT},
\[ \sup_{0\le \theta\le 2\pi}{\|x+e^{i\theta} y\|} \ge \frac{1}{2\pi}\int_0^{2\pi}\|x+e^{i\theta} y\| d\theta
\ge \left\| \left(|x|^2 + \frac12|y|^2\right)^{1/2}\right\|. \] By
Lemma~2.3 in \cite{Lee}, for every nonzero pair $x,y$ in $X$, there
exist $\delta_1=\delta_1(\|x\|, \|y\|)>0$ and  $z\in X$ with $0\le
z\le |y|$ and $\|z\|\ge  \|y\|/2$ such that the following holds:
\[ \| (|x|^2 + |y|^2)^{1/2} \| \ge \|\ |x| + \delta_1 |z|\ \|.\]
Hence for every $x\in S_X$ and $\epsilon>0$, we get
\[ \inf \left\{ \sup_{0\le \theta\le 2\pi}{\|x+e^{i\theta} y\|}: \|y\|\ge \epsilon \right\}
\ge \inf \left\{ \|\ |x| + |y|\ \|: \|y\|\ge \frac
\epsilon{2\sqrt2}\delta_1(1,\frac\epsilon{\sqrt{2}}) \right\}.\]
Hence the upper local uniform monotonicity implies the local
c-convexity. The proof is complete.
\end{proof}

A function $\varphi: \mathbb{R}\rightarrow [0,\infty]$ is said to be
an {\it Orlicz function} if $\varphi$ is even, convex continuous and
vanishing only at zero. Let $w=\{w(n)\}$ be a {\it weight sequence}
,that is, a non-increasing sequence of positive real numbers
satisfying $\sum_{n=1}^\infty w(n) =\infty$. Given a sequence $x$,
$x^*$ is the decreasing rearrangement of $|x|$.

The Orlicz-Lorentz sequence space $\lambda_{\varphi,w}$ consists of
all sequences $x=\{x(n)\}$ such that for some $\lambda>0$,
\[
\varrho_\varphi(\lambda x)= \sum_{n=1}^\infty \varphi(\lambda
x^*(n))w(n) < \infty,
\]
and equipped with the norm $\|x\| = \inf\{\lambda>0:
\varrho_\varphi(x/\lambda)\le 1\}$, $\lambda_{\varphi,w}$ is a
Banach space. We say an Orlicz function $\varphi$ satisfies the {\it
condition $\delta_2$} $(\varphi\in \delta_2)$ if there exist $K>0$,
$u_0>0$ such that $\varphi(u_0)>0$ and the inequality
\[ \varphi(2u) \le K \varphi(u)\] holds for $u\in [0, u_0]$.

It was proved in \cite[Corollary~4]{FK} that the Orlicz-Lorentz
sequence space $\lambda_{\varphi,w}$ is strictly monotone if and
only if it is both upper and lower locally uniformly monotone. They
also showed that the strict monotonicity of $\lambda_{\varphi,w}$ is
equivalent to the fact that $\varphi \in \delta_2.$ In this case,
the Orlicz-Lorentz sequence space $\lambda_{\varphi,w}$ is locally
c-convex by Proposition \ref{propmonotonecomplex} and order
continuous by Theorem 2 and Corollary 4 of \cite{FK}. If
$\varphi(u)=|u|^p$ for some $1\le p<\infty$ and if $w\equiv 1$, then
$\lambda_{\varphi, w} = \ell_p$. Hence we obtain the following
corollary by Theorem~\ref{thmsilov} which extends a result in
\cite{ACLQ}.

\begin{cor}\label{Orlicz}  Given  an Orlicz function $\varphi \in \delta_2$
and  a weight sequence $w,$
\[\partial \mathcal{A}_u(B_{\lambda_{\varphi, w}}) =
\partial \mathcal{A}_b (B_{\lambda_{\varphi, w}}) = S_{{\lambda_{\varphi, w}}}.\]
\end{cor}

\section{Boundaries for $\mathcal{A}_b(B_X)$}

Recall that a Banach sequence space $X$ is called {\it rearrangement
invariant} if $y\in X$ and $\|y\|=\|x\|$ whenever $y$ is a sequence
with $y^*=x^*$ for some $x\in X.$ Let $X$ be a rearrangement
invariant Banach sequence space. Given any finite subset $M$ of
natural numbers, let $\phi:\mathbb{N}\rightarrow \mathbb{N}\setminus
M$ be the order preserving bijection and let $P_M$ be the isometry
from $\{x\in X: {\rm supp}( x) \cap M =\emptyset\}$ onto $X$ given
by
\[ P_M(x) = \sum_{i=1}^\infty \langle x, e_{\phi(i)}\rangle
e_i,\] where the sum is a formal series and  $\mathrm{supp}(x)=\{k
\in \mathbb{N} ~:~ x(k) \neq 0\}.$ If $\mathrm{supp}(x)$ is finite,
$x$ is called a {\it finite vector}. Now assume that $X$ has the
following additional property:
\\ For each finite vector $x\in B_X,$
there exist $\epsilon=\epsilon(x)>0$ such that for all $y \in B_X$
with $\mathrm{supp}(x) \cap \mathrm{supp}(y)=\emptyset$,
\begin{equation}\label{eqspecicalcase} \| x + \epsilon y \| \le
1.\end{equation} For each finite vector $x\in S_X,$ let $\eta(x)>0$
be the supremum of the set of all $\epsilon>0$ satisfying
(\ref{eqspecicalcase}). Observe that $\eta(x) \le 1$ and
\[ \mathcal{F}(x) \supset\{ x+ \eta(x)y :  y \in B_X,
\ {\rm supp}(x) \cap {\rm supp} (y) =\emptyset\}.\]

If $X$ is a rearrangement invariant Banach sequence space satisfying
the property (\ref{eqspecicalcase}) and if $x \in
Ext_\mathbb{C}(B_F)$, where $F$ is the subspace spanned by a finite
number of $\{e_k\}_k$, then for all $y \in X$ with $x+y\in
\mathcal{F}(x)$, we have
\[{\rm supp}(x) \cap {\rm supp}(y)=\emptyset,\  \mbox{and}\ \  \|x+\eta(x)y\|\le 1.\]

Let $S\subset B_X$ and let $F$ be the subspace spanned by a finite
number of $\{e_k\}_k.$ Given $x \in B_F$ we define
\[ S(x) = \{ P_{{\mathrm{supp}(x)}}(y):  y\in B_X, ~~x+ \eta(x)y \in S,
 ~~\mathrm{supp}(x) \cap \mathrm{supp}(y)
=\emptyset\} .\] Put, for each $0<\epsilon<1$,
\[ C(S, \epsilon) = \sup\{ |f(0)| : f\in \mathcal{A}_b(B_X),\  \|f\|\le
1,\ |f(z)|<1-\epsilon \mathrm{~~for~~all~~} z\in S\}.\] Then $S$ is
called a {\it 0-boundary} for $\mathcal{A}_b (B_X)$ if
$C(S,\epsilon)<1$ for every $\epsilon>0$. A family
$\{S_\gamma\}_{\gamma\in \Gamma}$ of subsets of $B_X$ is called a
{\it uniform family of 0-boundaries} for $\mathcal{A}_b (B_X)$ if
$\sup_{\gamma\in \Gamma} C(S_\gamma, \epsilon)<1$ for every
$\epsilon>0$.

\begin{thm} \label{Ab1-boundary}
Let $X$ be a rearrangement invariant Banach sequence space
satisfying property (\ref{eqspecicalcase}) and let $V$ be a boundary
for $\mathcal{A}_b(B_X)$ consisting of norm-one finite vectors.
Assume also that $S \subset B_X$ has the property that
$\{S(x)\}_{x\in V}$ is a uniform family of 0-boundaries for
$\mathcal{A}_b(B_X).$ Then $S$ is a boundary for
$\mathcal{A}_b(B_X).$
\end{thm}
\begin{proof}
Suppose $S$ is not a boundary for $\mathcal{A}_b(B_X).$ Then there
is $f\in \mathcal{A}_b(B_X)$ with $\|f\|=1$ and $0 < \delta < 1$
such that $|f(z)|<1-\delta$ for all $z \in S$. The assumption on $V$
implies that there exists a sequence $\{x_n\}_{n=1}^{\infty}$ in $V$
such that
$$\lim_{n \rightarrow \infty} |f(x_n)|=1.$$ For each $n \in
\mathbb{N},$ $\mathrm{supp}(x_n)$ is finite and there exists
$\eta(x_n)>0$ such that for every $y \in B_X$ with
$\mathrm{supp}(x_n) \cap \mathrm{supp}(y)=\emptyset$ $$ \| x_n+
\eta(x_n) y\| \le 1.$$ Define $\phi_n$ on $B_X$
\[ x \mapsto f \left( x_n + \eta(x_n) P_{\mathrm{supp}(x_n)}^{-1}(x) \right),\]
Then $\phi_n \in \mathcal{A}_b(B_X)$ and $\|\phi_n\|\le 1$ for all
$n.$ Moreover, for each $x \in S(x_n)$
\[ |\phi_n(x)|<1-\delta \ \ \mbox{and}\ \
\lim_{n\rightarrow\infty} |\phi_n(0)| = 1,\] and this contradicts
the assumption that $\{S(x)\}_{x \in V}$ is a uniform family of
0-boundaries for $\mathcal{A}_b(B_X)$.
\end{proof}

We shall use the following two lemmas which are proved in \cite{G1}.
\begin{lem}\cite{G1}\label{lemG1}
Let $0\le r<1$ and assume $\{S_\gamma\}_{\gamma\in \Gamma}$ is a
family of subsets of $B_X$ such that $S_\gamma \cap rB_X\neq
\emptyset$ for each $\gamma\in \Gamma$. Then
$\{S_\gamma\}_{\gamma\in \Gamma}$ is a uniform family of
0-boundaries for $\mathcal{A}_b(B_X).$
\end{lem}

\begin{lem}\cite{G1}\label{lemG2}
Let $\theta_0>0$ and let $\{S_\gamma\}_{\gamma\in \Gamma}$ be a
family of subsets of $B_X$ with the following property: for each
$\gamma\in \Gamma$ there is some $x_\gamma \in S_\gamma$ such that
$e^{i\theta}x_\gamma\in S_\gamma$ for every $|\theta|\le \theta_0$.
Then $\{S_\gamma\}_{\gamma\in \Gamma}$ is a uniform family of
0-boundaries for $\mathcal{A}_b(B_X).$
\end{lem}

Theorem \ref{Ab1-boundary} and Lemma \ref{lemG1} show the following
corollary.
\begin{cor}Let $X$ be a rearrangement invariant Banach
sequence space satisfying property (\ref{eqspecicalcase}) and let
$V$ be a boundary for $\mathcal{A}_b(B_X)$ consisting of norm-one
finite vectors. Assume that $S \subset B_X$ and  that there is $0
\le r < 1$ such that for each $x \in V$ there exists $y \in X$ such
that
$$ \|y\| \le r ~,~ \mathrm{supp}(x) \cap \mathrm{supp}(y)=\emptyset \mathrm{~~and~~}
x+\eta(x)y \in S.$$ Then $S$ is a boundary for $\mathcal{A}_b(B_X).$
\end{cor}

By Theorem \ref{Ab1-boundary} and Lemma \ref{lemG2}, we get the
following corollary.
\begin{cor}Let $X$ be a rearrangement invariant Banach
sequence space satisfying property (\ref{eqspecicalcase}) and let
$V$ be a boundary for $\mathcal{A}_b(B_X)$ consisting of norm-one
finite vectors. Assume that $S \subset B_X$ and assume that there is
$\theta_0 >0$ such that for each $x \in V$ there exists $y \in B_X$
such that
$$\mathrm{supp}(x) \cap \mathrm{supp}(y)=\emptyset \mathrm{~~and~~}
x+\eta(x) e^{i \theta} y \in S \mathrm{~~for~~all~~} |\theta| \le
\theta_0.$$ Then $S$ is a boundary for $\mathcal{A}_b(B_X).$
\end{cor}

\begin{example}
Assume that $\psi=\{\psi(n)\}$ is an strictly increasing sequence
with $\psi(0)=0$, $\psi(n)
> 0$ for $n\in \mathbb N$. The {\it Marcinkiewicz sequence space}
$m_\psi$ consists of all sequences $x= \{x(n)\}$ such that
\[ \|x\|_{m_\psi}=\sup_{n \in \mathbb{N}}\frac{\sum_{k=1}^{n}
x^*(k)}{\psi(n)}< \infty.\] Let $m_\Psi^0$ be the closed subspace of
$m_\psi$, equipped with the same norm $\|\cdot\|_{m_\psi}$
consisting of all $x \in m_\psi$ satisfying
\[ \lim_{n \rightarrow \infty} \frac{\sum_{k=1}^{n}
x^*(k)}{\psi(n)} = 0.\] Without loss of generality we can add (and
we will) in the above definition the assumption that the sequence
$\{ \frac{\psi(n)} {n}\}_{n=1}^{\infty}$ is decreasing \cite{KL}.
Notice that if $\psi(n)=n$, then $m_\psi=\ell_\infty$ and
$m_\psi^0=c_0$, and if $\lim_n\psi(n)<\infty$, then
$m_\psi^0=\{0\}$.

It is shown in \cite{KL1} that if $\lim_n \psi(n)=\infty$, then for
each $x\in B_{m_\psi^0}$, there exist $n \in \mathbb{N}$ and
$\epsilon
>0$ such that $\| x+\lambda y\| \le 1$ for all $y \in B_{m_\psi}$
with $y=(0,\cdots,0,y(n+1),y(n+2),\cdots)$ and all $\lambda$ with
$|\lambda| \le \epsilon$. Now it is easy to see that $m_\psi^0$
satisfies (\ref{eqspecicalcase}) because $m_\psi^0$ is a
rearrangement invariant sequence space.
\end{example}

\subsection*{Acknowledgments} The authors thank the
referee whose careful reading and suggestions led to a much
improved version of this paper.

\bibliographystyle{amsplain}

\end{document}